\documentclass[12pt]{amsart}
\usepackage{amssymb}
\begin{document}

\def\dbl{[\hskip -1pt[}
\def\dbr{]\hskip -1pt]}

\title[Convergence of formal embeddings in codimension one]{Convergence of formal embeddings
between\\ real-analytic hypersurfaces in codimension one}
\author{Nordine Mir}
\address{Universit\'e de Rouen, Laboratoire de Math\'ematiques Rapha\"el Salem, UMR 6085 CNRS, 76821 Mont-Saint-Aignan Cedex, France}
\email{Nordine.Mir@univ-rouen.fr}
\thanks{\noindent 2000 {{\em Mathematics Subject Classification.}  32H02, 32V20, 32V30,
32V40}}

\def\Label#1{\label{#1}}
\def\1#1{\ov{#1}}
\def\2#1{\widetilde{#1}}
\def\6#1{\mathcal{#1}}
\def\4#1{\mathbb{#1}}
\def\3#1{\widehat{#1}}
\def\K{{\4K}}
\def\LL{{\4L}}

\def\C{{\4C}}
\def\R{{\4R}}
\def \MM{{\4M}}

\def\Re{{\sf Re}\,}
\def\Im{{\sf Im}\,}

\numberwithin{equation}{section}
\def\s{s}
\def\k{\kappa}
\def\ov{\overline}
\def\span{\text{\rm span}}
\def\ad{\text{\rm ad }}
\def\tr{\text{\rm tr}}
\def\xo {{x_0}}
\def\Rk{\text{\rm Rk\,}}
\def\sg{\sigma}
\def \emxy{E_{(M,M')}(X,Y)}
\def \semxy{\scrE_{(M,M')}(X,Y)}
\def \jkxy {J^k(X,Y)}
\def \gkxy {G^k(X,Y)}
\def \exy {E(X,Y)}
\def \sexy{\scrE(X,Y)}
\def \hn {holomorphically nondegenerate}
\def\hyp{hypersurface}
\def\prt#1{{\partial \over\partial #1}}
\def\det{{\text{\rm det}}}
\def\wob{{w\over B(z)}}
\def\co{\chi_1}
\def\po{p_0}
\def\fb {\bar f}
\def\gb {\bar g}
\def\Fb {\ov F}
\def\Gb {\ov G}
\def\Hb {\ov H}
\def\zb {\bar z}
\def\wb {\bar w}
\def \qb {\bar Q}
\def \t {\tau}
\def\z{\chi}
\def\w{\tau}
\def\Z{\zeta}
\def\phi{\varphi}
\def\eps{\varepsilon}

\def \T {\theta}
\def \Th {\Theta}
\def \L {\Lambda}
\def\b {\beta}
\def\a {\alpha}
\def\o {\omega}
\def\l {\lambda}

\def \im{\text{\rm Im }}
\def \re{\text{\rm Re }}
\def \Char{\text{\rm Char }}
\def \supp{\text{\rm supp }}
\def \codim{\text{\rm codim }}
\def \Ht{\text{\rm ht }}
\def \Dt{\text{\rm dt }}
\def \hO{\widehat{\mathcal O}}
\def \cl{\text{\rm cl }}
\def \bR{\mathbb R}
\def \bS{\mathbb S}
\def \bK{\mathbb K}
\def \bD{\mathbb D}
\def \bC{\mathbb C}
\def \C{\mathbb C}
\def \N{\mathbb N}
\def \bL{\mathbb L}
\def \bZ{\mathbb Z}
\def \bN{\mathbb N}
\def \scrF{\mathcal F}
\def \scrK{\mathcal K}
\def \mc #1 {\mathcal {#1}}
\def \scrM{\mathcal M}
\def \cR{\mathcal R}
\def \scrJ{\mathcal J}
\def \scrA{\mathcal A}
\def \scrO{\mathcal O}
\def \scrV{\mathcal V}
\def \scrL{\mathcal L}
\def \scrE{\mathcal E}
\def \hol{\text{\rm hol}}
\def \aut{\text{\rm aut}}
\def \Aut{\text{\rm Aut}}
\def \J{\text{\rm Jac}}
\def\jet#1#2{J^{#1}_{#2}}
\def\gp#1{G^{#1}}
\def\gpo{\gp {2k_0}_0}
\def\emmp {\scrF(M,p;M',p')}
\def\rk{\text{\rm rk\,}}
\def\Orb{\text{\rm Orb\,}}
\def\Exp{\text{\rm Exp\,}}
\def\Span{\text{\rm span\,}}
\def\d{\partial}
\def\D{\3J}
\def\pr{{\rm pr}}

\def \CZZ {\C \dbl Z,\zeta \dbr}
\def \D{\text{\rm Der}\,}
\def \Rk{\text{\rm Rk}\,}
\def \CR{\text{\rm CR}}
\def \ima{\text{\rm im}\,}
\def \I {\mathcal I}

\def \M {\mathcal M}

\newtheorem{Thm}{Theorem}[section]
\newtheorem{Cor}[Thm]{Corollary}
\newtheorem{Pro}[Thm]{Proposition}
\newtheorem{Lem}[Thm]{Lemma}

\theoremstyle{definition}\newtheorem{Def}[Thm]{Definition}

\theoremstyle{remark}
\newtheorem{Rem}[Thm]{Remark}
\newtheorem{Exa}[Thm]{Example}
\newtheorem{Exs}[Thm]{Examples}

\def\bl{\begin{Lem}}
\def\el{\end{Lem}}
\def\bp{\begin{Pro}}
\def\ep{\end{Pro}}
\def\bt{\begin{Thm}}
\def\et{\end{Thm}}
\def\bc{\begin{Cor}}
\def\ec{\end{Cor}}
\def\bd{\begin{Def}}
\def\ed{\end{Def}}
\def\br{\begin{Rem}}
\def\er{\end{Rem}}
\def\be{\begin{Exa}}
\def\ee{\end{Exa}}
\def\bpf{\begin{proof}}
\def\epf{\end{proof}}
\def\ben{\begin{enumerate}}
\def\een{\end{enumerate}}

\keywords{}

\maketitle



\section{Introduction and results}\Label{int}

By a classical result of Chern-Moser \cite{CM}, any formal
biholomorphic transformation in the complex $N$-dimensional space,
$N\geq 2$, sending two real-analytic strongly pseudoconvex
hypersurfaces into each other is in fact convergent, i.e.\ given
by the power series of a local holomorphic map. Several
generalizations of this result have recently been established for
more general classes of real-analytic hypersurfaces in the
equidimensional case (see e.g.\ \cite{BERsurvey, M00} and the
references therein). On the other hand, it is conjectured that
Chern-Moser's result can be extended to formal embeddings sending
a real-analytic strongly pseudoconvex hypersurface $M\subset \C^N$
into another such hypersurface $M'\subset \C^{N'}$ with $N'>N$.
Here we recall that by a formal embedding $F\colon (\C^{N},p)\to
(\C^{N'},p')$ sending $M$ into $M'$, $p\in M$, $p'\in M'$, we mean
a formal holomorphic map for which the pullback under $F$ of any
local real-analytic defining function for $M'$ near $p'$ vanishes
on $M$ (as a formal power series), and for which the induced
differential $dF(p)|_{\C T_pM}\colon \C T_pM\to \C T_{p'}M'$ is
injective. The main difficulty that one has to face in such a
setting, and that was pointed out by a number of authors (see
e.g.\ \cite{F, HH}), lies in the codimension $N'-N>0$ of the
embedding and indeed, up to now, only partial results on the above
question have been obtained under various additional assumptions
on the mapping or manifolds (see e.g.\ \cite{la, MMZ3}). For
instance, it follows from the recent results of \cite{MMZ3} that
the above conjecture holds if the target hypersurface is moreover
assumed to be real-algebraic. In this paper, we make a step toward
the understanding of the above problem by giving a complete
solution in the one-codimensional case. Indeed, we have:

\begin{Thm}\Label{main1} Any formal embedding sending a
real-analytic strongly pseudoconvex hypersurface $M\subset \C^N$
into another such hypersurface $M'\subset \C^{N+1}$ is convergent.
\end{Thm}

Theorem \ref{main1} gives also the first positive answer to a long
standing open problem which consists in providing a regularity
result for embeddings of positive codimension between
real-analytic strongly pseudoconvex hypersurfaces that holds at
{\em all} points of the source manifold (see e.g.\ \cite{F, HH}).
As a byproduct of the proof of Theorem \ref{main1}, we also obtain
two independent convergence results for formal maps of
Levi-nondegenerate real-analytic hypersurfaces. In what follows,
we consider  formal maps $F\colon (\C^N,p)\to (\C^{N'},p')$
sending $M$ into $M'$ that are CR {\em transversal} i.e.\ for
which $dF(p)(\C T_pM)\not \subset T^{1,0}_{p'}M'\oplus
T^{0,1}_{p'}M'$ where $T_{p'}^{1,0}M'\subset \C T_{p'}M'$ (resp.\
$T_{p'}^{0,1}M'\subset \C T_{p'}M'$) denotes the $(1,0)$ (resp.\
$(0,1)$) tangent space of $M'$ at $p'$, and we shall prove the
following.

\begin{Thm}\Label{main3}
Any formal CR transversal map sending a real-analytic
Levi-nondegenerate hypersurface $M\subset \C^N$ into another such
hypersurface $M'\subset \C^{N+1}$ is convergent.
\end{Thm}

As is well-known, any formal CR transversal map in the setting of
Theorem \ref{main3} is necessarily an embedding. Moreover, it
turns out that the condition of CR transversality is automatically
satisfied by all formal embeddings between real-analytic strongly
pseudoconvex hypersurfaces (see e.g.\ \cite{la,EL}) and therefore,
Theorem \ref{main1} follows immediately from Theorem \ref{main3}.
In addition, it is easy to see that Theorem \ref{main3} can not be
extended to formal embeddings of codimension higher or equal to
two (see e.g.\ \cite{la}). For those embeddings that are not
necessarily CR transversal, we also have:

\begin{Thm}\Label{main2}
Let  $M\subset \C^N$ and $M'\subset \C^{N+1}$ be real-analytic
Levi-nondegenerate hypersurfaces through points $p$ and $p'$
respectively. Assume that $M'$ does not contain any $($smooth$)$
complex curve through $p'$. Then any formal embedding $F\colon
(\C^N,p)\to (\C^{N+1},p')$ sending $M$ into $M'$ is convergent.
\end{Thm}

The proof of Theorems \ref{main3} and \ref{main2}, that will be
given in \S \ref{proof}, is inspired by the Segre variety approach
developped by Webster \cite{W} and the CR vector field approach
carried out in the works of Huang \cite{Hu94} and Baouendi,
Ebenfelt and Rothschild \cite{BERbook}. As a preliminary step, we
need to establish in \S \ref{crit} a useful criterion for the
convergence of a formal power series that satisfies a certain type
of identity (see Proposition \ref{trick}). The proof of such a
statement also uses some arguments from our previous works
\cite{M00, M02}.

\section{A criterion for the convergence of a formal power
series}\Label{crit}

 Throughout the paper,
we denote by $\C\dbl x\dbr$ (resp.\ $\C \{x\}$),
$x=(x_1,\ldots,x_k)$, the ring of formal (resp.\ convergent) power
series in $k$ indeterminates with complex coefficients. If $x^0\in
\C^k$, $\C\dbl x-x^0\dbr$ and $\C \{x-x^0\}$ will denote the
corresponding ring of series centered at $x^0$. Given a formal
power series $g(x)\in \C \dbl x\dbr$, we also denote by $\1g (x)$
the formal power series obtained from $g(x)$ by taking complex
conjugates of its coefficients. Given moreover a (germ at the
origin of a) complex submanifold $S\subset \C^k$, we write
$g(x)\equiv 0$ for $x\in S$ to mean that $g\circ \nu \equiv 0$ for
any parametrization $\nu$ of $S$. We start by stating the
following well-known lemma (see e.g.\ \cite[Propositions 4.2 and
6.2]{M00}).

\bl\Label{warm-up} Let $R(x,y)\in \C\{x,y\}$,
$x=(x_1,\ldots,x_k)$, $y\in \C$, $h(x)\in \C \dbl x\dbr$ with
$h(0)=0$ and $v\in (\C\{t\})^k$, $t=(t_1,\ldots,t_q)$. Then the
following holds:
\begin{enumerate}
\item[(i)] If $R(x,y)\not \equiv 0$ and $R(x,h(x))\equiv 0$ then $h(x)$ is convergent.
\item[(ii)] If $(h\circ v)(t)$ is convergent and $v$ is of generic
rank $k$, then $h(x)$ is itself convergent.
\end{enumerate}
 \el

Let $M\subset \C^N$ be a (germ of a) real-analytic hypersurface
through the origin, and $\rho (Z,\1Z)$ be a real-analytic defining
function for $M$ defined in a connected neighborhood $U$ of $0$ in
$\C^N$, with non-vanishing gradient on $U$. Recall that the
complexification $\M$ of $M$ is the complex submanifold of
$\C^{2N}$ defined as follows
\begin{equation}\Label{comp}
\M:=\{(Z,\zeta)\in U\times U^*: \rho (Z,\zeta)=0\},
\end{equation}
where for any subset $V\subset \C^N$, we have denoted
$V^*:=\{\1w:w\in V \}$. Recall also that $M$ is said to be of {\em
finite type} at the origin if there is no complex hypersurface
contained in $M$ through $0$ (see \cite{BERbook}). We may now
formulate one of the main tools in the proof of Theorems
\ref{main3} and \ref{main2}.

\begin{Pro}\Label{trick}
Let $M\subset \C^N$ be a real-analytic hypersurface of finite type
through the origin, and $\M\subset \C^N_Z\times \C^N_{\zeta}$ its
complexification as given by {\rm (\ref{comp})}. Let $H(Z)$ be a
formal power series, with $H(0)=0$, satisfying at least one of the
following conditions:
\begin{enumerate}
\item[(i)] there exists $G(\zeta):=(G_1(\zeta),\ldots,G_m(\zeta))$ a vector-valued
formal power series and $A(Z,\zeta,X,T)\in
\C\{Z,\zeta,X-G(0),T\}$, $X=(X_1,\ldots,X_m)$, $T\in \C$, such
that $A(Z,\zeta,G(\zeta),T)\not \equiv 0$ for $(Z,\zeta)\in \M$
and such that $A(Z,\zeta,G(\zeta),H(Z))\equiv 0$ for $(Z,\zeta)\in
\M$.
\item[(ii)]  there exists  $B(Z,\zeta,\2T,T)\in
\C\{Z,\zeta,\2T,T\}$, $T,\2T\in \C$, such that
$B(Z,\zeta,\2T,T)\not \equiv 0$ for $(Z,\zeta)\in \M$ and such
that $B(Z,\zeta,\1H(\zeta),H(Z))\equiv 0$ for $(Z,\zeta)\in \M$.
\end{enumerate}
Then the formal power series $H(Z)$ is necessarily convergent.
\end{Pro}

\br\Label{note}
 In the case where the power series $A,B$ in
Proposition \ref{trick} are polynomials in $T$ and $\2T$, the
above conclusion follows from \cite{M02, MMZ3}.\er

\begin{proof}
We choose local holomorphic coordinates $Z=(z,w)\in \C^{N-1}\times
\C$ so that $M$ is given near the origin by an equation of the
form
\begin{equation}\Label{hunt}
w=Q(z,\1z,\1w),
\end{equation}
for some holomorphic function $Q(z,\chi,\tau)$ defined in a
neighborhood of $0\in \C^{N-1}\times \C^{N-1}\times \C$. We may
also assume (see e.g.\ \cite{BERbook}) that $Q$ satisfies
\begin{equation}\Label{reality}
Q(z,\chi,\1Q(\chi,z,w))\equiv w,\ Q(z,0,w)=Q(0,\1z,w)\equiv 0.
\end{equation}
Consider the parametrizations $v^1,v^2,v^3$ of the Segre sets up
to order 3 attached to $M$ at the origin (see \cite{BERbook})
given by
\begin{equation}
\ v^{1}(t^1):=(t^1,0),\ v^2(t^1,t^2):=(t^2,Q(t^2,\1{v^1}(t^1))),\
v^3(t^1,t^2,t^3):=(t^3,Q(t^3,\1{v^2}(t^1,t^2))),
\end{equation}
where each $t^j\in \C^{N-1}$ is sufficiently close to the origin.
Recall that since $M$ is of finite type at $0$, the holomorphic
map
\begin{equation}\Label{eta}
\eta\colon (\C^{3N-3},0)\ni (t^1,t^2,t^3)\mapsto
(v^2(t^1,t^2),\1{v^3}(t^1,t^2,t^3))\in (\M,0)
\end{equation}
is of generic maximal rank $2N-1$ (see e.g.\ \cite{BERbook}).

Let $A,G,H$ be as in (i). Then after composing the given
identities with the map $(Z,\zeta)=\eta (t^1,t^2,t^3)$, we have
\begin{equation}\Label{sp}
A(v^2(t^1,t^2),\1{v^3}(t^1,t^2,t^3),(G\circ
\1{v^3})(t^1,t^2,t^3),(H\circ v^2)(t^1,t^2))\equiv 0,
\end{equation}
and
\begin{equation}\Label{st}
\Delta
(t^1,t^2,t^3,T):=A(v^2(t^1,t^2),\1{v^3}(t^1,t^2,t^3),(G\circ
\1{v^3})(t^1,t^2,t^3),T)\not \equiv 0,
\end{equation}
in view of the generic rank of $\eta$. From (\ref{st}) and
(\ref{sp}), we may choose a  multiindex $\beta_0\in \N^{N-1}$ such
that
\begin{equation}\Label{delta}
\2{\Delta}(t^1,t^2,T):=\left[\frac{\partial^{|\beta_0|}\Delta}{\partial
t_3^{\beta_0}}(t^1,t^2,t^3,T)\right]_{\big|{t^3=t^1}}\not \equiv
0\ {\rm and}\  \2{\Delta}(t^1,t^2,(H\circ v^2)(t^1,t^2))\equiv 0.
\end{equation}
Note that by using the identity $v^3(t^1,t^2,t^1)=v^1(t^1)$ (which
follows from (\ref{reality})), we may rewrite
$\2{\Delta}(t^1,t^2,T)=\widehat{\Delta}(t^1,t^2,
(((\partial^{\beta}G)\circ \1{v^1})(t^1))_{|\beta|\leq
|\beta_0|},T)$ for some holomorphic function $\widehat{\Delta}$ in
a neighborhood of $(0,0,(((\partial^{\beta}G)(0))_{|\beta|\leq
|\beta_0|},0)$. Now as in \cite{M00,M02}, by differentiating the
second identity in (\ref{delta}) with respect to $t^2$, setting
$t^2=0$, using the identity $v^2(t^1,0)=0$, and applying Artin's
approximation theorem \cite{A68}, we may find for any positive
integer $k$ a convergent power series mapping $Y^k(t^1)$ which
agrees up to order $k$ with $(((\partial^{\beta}G)\circ
\1{v^1})(t^1))_{|\beta|\leq |\beta_0|}$ at the origin and which
satisfies the identity
\begin{equation}\Label{fruit}
\widehat{\Delta}(t^1,t^2,Y^k(t^1),(H\circ v^2)(t^1,t^2))\equiv 0.
\end{equation}
Moreover, by choosing $k$ large enough, we may achieve the
condition $\widehat{\Delta}(t^1,t^2,Y^k(t^1),T)\not \equiv 0$ in
view of (\ref{delta}). We may therefore apply Lemma \ref{warm-up}
(i) to conclude that $H\circ v^2$ is convergent. Since $M$ is of
finite type at $0$, it is easy to see that the generic rank of
$v^2$ is $N$, and hence, from Lemma \ref{warm-up} (ii), it follows
that $H$ is convergent. The proof of Proposition \ref{trick} (i)
is complete.

Let $B,H$ as in (ii). Expand $B$ as a Taylor series as follows
$B(Z,\zeta,\2T,T)=\sum_{j=0}^{\infty}b_j(Z,\zeta,\2T)T^j$. There
are two cases to consider.

{\sc First case.} There exists $j_0$ such that
$b_{j_0}(Z,\zeta,\1H(\zeta))\not \equiv 0$ for $(Z,\zeta)\in \M$.
Then $B(Z,\zeta,\1H(\zeta),T)$ is nontrivial for $(Z,\zeta)\in \M$
and the convergence of $H$ follows from Proposition \ref{trick}
(i) proved above.

{\sc Second case.} For all $j$, $b_j(Z,\zeta,\1H(\zeta)) \equiv 0$
for $(Z,\zeta)\in \M$. Since $B(Z,\zeta,\2T,T)\not \equiv 0$ for
$(Z,\zeta)\in \M$, there exists $j_1$ such that
$B_{j_1}(Z,\zeta,\2T)\not \equiv 0$ for $(Z,\zeta)\in \M$.  The
convergence of the series $H$ then follows from Lemma
\ref{warm-up} (i) by noticing that
$\1{B_{j_1}}(\zeta,Z,H(Z))\equiv 0$ for $(Z,\zeta)\in \M$. This
completes the proof of Proposition \ref{trick} (ii).
\end{proof}

\section{Proofs of Theorems \ref{main2} and \ref{main3}}\Label{proof}

\begin{proof}[Proof of Theorem {\rm \ref{main2}}] Without loss of generality, we may assume that
$p=p'=0$. We choose local holomorphic coordinates $Z'=(z',w')\in
\C^{N}\times \C$ near the origin so that the Levi-nondegenerate
real-analytic hypersurface $M'$ is given by the equation
\begin{equation}\Label{M'}
\1{w'}=\Theta(\1{z'},{z'},{w'}),
\end{equation}
where $\Theta=\Theta(\chi',Z')$ is a holomorphic function near the
origin in $\C^{N}\times \C^{N+1}$. We set
$\rho'(Z',\1{Z'}):=\1{w'}-\Theta(\1{z'},{z'},{w'})$ and may assume
(see \cite{CM}) that $\Theta(\chi',z',0)$ vanishes at the origin
up to order one and also satisfies
\begin{equation}\Label{normalisation}
\Theta_{\chi'_jz_{k}'}(0)=0,\quad {\rm for}\ j\not =k,\quad
\Theta_{\chi'_jz_{j}'}(0)=\pm 1,\quad j,k=1,\ldots,N.
\end{equation}
 In these coordinates, we split the formal map $F$ as follows
$F=(f,g)\in \C^{N}\times \C$, where $f=(f_1,\ldots,f_{N})$. At the
source $\C^N$ space, we denote our coordinates by $Z=(z,w)\in
\C^{N-1}\times \C$ which we may assume to be normal coordinates as
in (\ref{hunt}) and (\ref{reality}). Since $F$ sends $M$ into
$M'$, we have the formal identity
\begin{equation}\Label{fundamental}
\1{g}(\zeta)=\Theta(\1{f}(\zeta),F(Z)),\ {\rm for}\ (Z,\zeta)\in
\M,
\end{equation}
where $\M$ is the complexification of $M$ as given by
(\ref{comp}). Let $\6L_1,\ldots,\6L_{N-1}$ be a basis of
holomorphic vector fields (in a neighborhood of $0$ in $\C^N\times
\C^N$, with holomorphic coefficients in $(Z,\zeta)$) tangent to (a
neighborhood of $0$ in) $\M$ that annihilate the projection
$\C^N\times \C^N \ni (Z,\zeta)\mapsto Z\in \C^N$. Applying each
$\6L_j$ to (\ref{fundamental}), we obtain
\begin{equation}\Label{differ}
\6L_{j}\1{g}(\zeta)=\sum_{k=1}^{N}\Theta_{\chi_{k}'}(\1f(\zeta),F(Z))\6L_j\1f_{k}(\zeta).
\end{equation}

By (\ref{differ}) and our choice of $\Theta$, we have
$(\6L_j\1g)(0)=0$ for all $j$ and therefore,  since $F$ is an
embedding, the rank of the matrix
$\left(\6L_j\1f_{k}(\zeta)\right)_{1\leq j\leq N-1 \atop 1\leq
k\leq N}$ is equal to $N-1$ at the origin. Then after
interchanging the components of $f$ if necessary, we may assume
that the rank of the matrix
$\left(\6L_j\1f_{k}(\zeta)\right)_{1\leq j\leq N-1\atop 1\leq
k\leq N-1}$ equals $N-1$ at the origin. Therefore, by using
Cramer's rule to (\ref{differ}), taking the complex conjugate of
(\ref{fundamental}), we obtain the system of formal equations

\begin{equation}\Label{hurry-up}
\left\{\aligned
&g(Z)=\1{\Theta}(f(Z),\1F(\zeta))\\
&\Theta_{\chi_{k}'}(\1f(\zeta),F(Z))=\Theta_{\chi_{N}'}(\1f(\zeta),F(Z))\,
P_k((\6L_j\1F(\zeta))_{1\leq j\leq
N-1})+S_k((\6L_j\1F(\zeta))_{1\leq j\leq N-1}),
\endaligned\right.
\end{equation}
for $k=1,\ldots,N-1$ and $(Z,\zeta)\in \M$. Here each $P_k,S_k$ is
a convergent power series centered at $((\6L_j\1F)(0))_{1\leq
j\leq N-1}$ (that is even rational). In view of
(\ref{normalisation}), we may solve the system (\ref{hurry-up}) by
making use of the implicit function theorem to obtain the
vectorial formal identity
\begin{equation}\Label{loi}
(f_1(Z),\ldots,f_{N-1}(Z),g(Z))=\Psi(Z,\zeta,((\partial^{\alpha}\1F)(\zeta))_{|\alpha|\leq
1},f_{N}(Z)),\ {\rm for}\ (Z,\zeta)\in \M,
\end{equation}
where
$\Psi=(\Psi_1,\ldots,\Psi_{N-1},\Psi_{N+1})=\Psi(Z,\zeta,\Lambda,T)\in
\left(\C\{Z,\zeta,\Lambda-((\partial^{\alpha}\1F)(0))_{|\alpha|\leq
1},T\}\right)^{N}$. (Here we have used the fact that the vector
fields $\6L_j$, $j=1,\ldots,N-1$, have holomorphic coefficients in
$(Z,\zeta)$.) For $k\in \{1,\ldots,N+1\}$, $k\not =N$, we write
the Taylor expansion
\begin{equation}\Label{taylor}
\Psi_k(Z,\zeta,\Lambda,T)=\sum_{i=0}^{\infty}\varphi_{k,i}
(Z,\zeta,\Lambda)T^i.
\end{equation}
 Applying
again each vector field $\6L_j$ to (\ref{loi}), we obtain for all
$k$ as above

\begin{equation}\Label{funny}
\begin{aligned}
0&=\6L_j\left(\Psi_k(Z,\zeta,((\partial^{\alpha}\1F)(\zeta))_{|\alpha|\leq
1},f_{N}(Z))\right)\\
&= \sum_{i=0}^{\infty}\6L_j \left(\varphi_{k,i}
(Z,\zeta,((\partial^{\alpha}\1F)(\zeta))_{|\alpha|\leq
1})\right)\left(f_{N}(Z)\right)^i,\ (Z,\zeta)\in \M.
\end{aligned}
\end{equation}
Note that there is a convergent power series
$\2{\Psi}_{k,j}(Z,\zeta,\widehat \Lambda,T)\in
\C\{Z,\zeta,\widehat\Lambda-((\partial^{\alpha}\1F)(0))_{|\alpha|\leq
2},T\}$ such that
\begin{equation}\Label{crap}
\2{\Psi}_{k,j}(Z,\zeta,((\partial^{\alpha}\1F)(\zeta))_{|\alpha|\leq
2},f_{N}(Z))=\6L_j\left(\Psi_k(Z,\zeta,((\partial^{\alpha}\1F)(\zeta))_{|\alpha|\leq
1},f_{N}(Z))\right).
\end{equation}
 Now, as in \cite{Hu94}, we come to a
dichotomy which will give the convergence of the map $F$ (in any
case).

{\sc First case.} There exist indices $k_0,j_0$ and $i_0$ such
that $\6L_{j_0} \left(\varphi_{k_0,i_0}
(Z,\zeta,((\partial^{\alpha}\1F)(\zeta))_{|\alpha|\leq
1}\right)\not \equiv 0$ for $(Z,\zeta)\in \M$. Then the formal
power series
$\2{\Psi}_{k_0,j_0}(Z,\zeta,((\partial^{\alpha}\1F)(\zeta))_{|\alpha|\leq
2},T)$ is nontrivial for $(Z,\zeta)\in \M$ and satisfies
$\2{\Psi}_{k_0,j_0}(Z,\zeta,((\partial^{\alpha}\1F)(\zeta))_{|\alpha|\leq
2},f_{N}(Z))\equiv 0$ on $\M$. Since $M$ is Levi-nondegenerate and
therefore of finite type at $0$, the convergence of $f_{N}$ then
follows from Proposition \ref{trick} (i). Using (\ref{loi}), the
now established convergence of $f_{N}$ and Proposition \ref{trick}
(i), we obtain the convergence of all other components of the map
$F$.

{\sc Second case.} For all indices $k,j$ and $i$ we have $\6L_{j}
\left(\varphi_{k,i}
(Z,\zeta,((\partial^{\alpha}\1F)(\zeta))_{|\alpha|\leq 1}\right)
\equiv 0$ for $(Z,\zeta)\in \M$. This means that for any $j\in
\{1,\ldots,N-1\}$, $\6L_j\left(\Psi
(Z,\zeta,((\partial^{\alpha}\1F)(\zeta))_{|\alpha|\leq
1},T)\right)\equiv 0$ for $(Z,\zeta)\in \M$. From this, it is then
easy to see that there exists a $\C^{N}$-valued formal power
series mapping $\Phi (Z,T)$ such that
\begin{equation}\Label{yes}
\Phi (Z,T)=\Psi
(Z,\zeta,((\partial^{\alpha}\1F)(\zeta))_{|\alpha|\leq 1},T),\
(Z,\zeta)\in \M.
\end{equation}
By (\ref{yes}) and Proposition \ref{trick} (i) applied to the
real-analytic hypersurface $M\times \C\subset \C^{N+1}$ (that is
of finite type), we obtain that $\Phi (Z,T)$ defines a holomorphic
map near $0\in \C^{N+1}$. We may therefore rewrite (\ref{loi}) as
follows
\begin{equation}\Label{loibis}
(f_1(Z),\ldots,f_{N-1}(Z),g(Z))=\Phi(Z,f_{N}(Z)).
\end{equation}
For $(Z,T)\in \C^N\times \C$ sufficiently close to the origin, we
set
\begin{equation}\Label{fest}
h(Z,T):=(\Phi_1(Z,T),\ldots,\Phi_{N-1}(Z,T),T,\Phi_{N+1}(Z,T))\in
\C^{N+1}. \end{equation}
 Consider the holomorphic function
defined near the origin in $\C^{N}\times \C^N\times \C\times \C$
by
\begin{equation}\Label{R}
R(Z,\zeta,T,\2T):=\rho'(h(Z,T),\1h(\zeta,\2T)).
\end{equation}
Note that necessarily $R(Z,\zeta,T,\2T)$ does not vanish
identically for $(Z,\zeta,T,\2T)\in (\M\times \C^2,0)$ since
otherwise $M'$ would contain a complex curve through the origin
which is impossible by assumption. In view of (\ref{fest}) and
(\ref{loibis}) we have $h(Z,f_N(Z))=F(Z)$ and therefore since $F$
sends $M$ into $M'$, we
 have
\begin{equation}\Label{ouf}
R(Z,\zeta,f_{N}(Z),\1{f_{N}}(\zeta))\equiv 0,\ {\rm for}\
(Z,\zeta)\in \M.
\end{equation}
Then by applying Proposition \ref{trick} (ii) to (\ref{ouf}), we
obtain that $f_{N}$ is convergent and hence $F$ too in view of
(\ref{loibis}). The proof of Theorem \ref{main2} is complete.
\end{proof}

\begin{proof}[Proof of Theorem {\rm \ref{main3}}] Let $F\colon (\C^N,p)\to
(\C^{N+1},p')$ be a given formal CR transversal map sending $M$
into $M'$. As in the proof of Theorem \ref{main2}, we may assume
that $p=p'=0$. Since $M$ and $M'$ are Levi-nondegenerate and $F$
is CR transversal, $F$ is a formal embedding (see e.g.\
\cite{la}). Therefore we may follow the same proof as that of
Theorem \ref{main2} and we notice that the only place (in that
proof) where the fact that the target hypersurface does not
contain any complex curve through $0$ is used is to show that the
holomorphic function $R$ given by (\ref{R}) does not vanish
identically for $(Z,\zeta,T,\2T)\in (\M \times \C^2,0)$. We now
show that this latter condition holds automatically under the
assumptions of Theorem \ref{main3}, which will finish its proof.
By contradiction suppose that $R(Z,\zeta,T,\2T)\equiv 0$ for
$(Z,\zeta,T,\2T)\in (\M\times \C^2,0)$. This means in view of
(\ref{R}) that $h(Z,T)\in M'$ for all $(Z,T)\in M \times \C$
sufficiently close to the origin. We now claim that $h$ is local
biholomorphism of $\C^{N+1}$. Indeed, by using e.g.\ the arguments
of \cite[Lemma 5.1]{EHZ}, we may assume, after composing the map
$F$ with an automorphism of the hyperquadric tangent to $M'$, that
the following normalization condition holds for the $N$-the
component of $F$:
\begin{equation}\Label{ebenfelt}
\frac{\partial F_N}{\partial z}(0)=0.
\end{equation}
Since $F$ is CR transversal we have $\partial g/\partial w(0)\not
=0$ (see e.g.\ \cite{la, EHZ}) and therefore, since $F$ is an
embedding, the map $Z\mapsto (f_1(Z),\ldots,f_{N-1}(Z),g(Z))$ is a
formal biholomorphism of $\C^N$. From this fact and the identities
(\ref{loibis}), (\ref{fest}) and (\ref{ebenfelt}), it is easy to
see that the Jacobian matrix of $h$ has indeed rank $N+1$ at the
origin, which proves the claim. We therefore have a local
biholomorphism $h$ of $\C^{N+1}$ satisfying $h(M\times \C)=M'$ (as
germs through the origin), which is impossible in view of the
Levi-nondegeneracy assumption on $M'$. The proof of Theorem
\ref{main3} is therefore complete.
\end{proof}

{\bf Acknowledgements.} The author would like to thank P.\
Ebenfelt, X.\ Huang and D.\ Zaitsev for helpful discussions and
many comments on a preliminary version of this paper. Special
thanks are due to P.\ Ebenfelt who pointed out an incomplete
statement in the original proof of Theorem \ref{main3} and
provided the correct argument.


\begin{thebibliography}{MMZ02a}


\bibitem{A68} {\sc Artin, M.} --- On the solutions of analytic equations. {\em Invent. math.}
{\bf 5} (1968), 277--291.

\bibitem{BERbook} {\sc Baouendi,~M.S.; Ebenfelt,~P.; Rothschild,~L.P.} --- {\em Real
Submanifolds in Complex Space and Their Mappings}. Princeton Math.
Series {\bf 47}, Princeton Univ. Press, 1999.

\bibitem{BERsurvey}{\sc Baouendi,~M.S.; Ebenfelt,~P.; Rothschild,~L.P.} --- Local geometric properties of
real submanifolds in complex space, {\em Bull. Amer. Math. Soc.},
{\bf 37}, (2000), 309--336.


\bibitem{CM} {\sc Chern,~S.S.; Moser,~J.K.} --- Real hypersurfaces in complex manifolds. {\em Acta Math.} {\bf 133} (1974), 219--271.

\bibitem{EL} {\sc Ebenfelt,~P.; Lamel,~B.} --- Finite jet
determination of constantly degenerate CR embeddings, (preprint
2002).


\bibitem{EHZ}{\sc Ebenfelt, P.; Huang, X.; Zaitsev, D.} --- The
equivalence problem and rigidity for hypersurfaces embedded into
hyperquadrics, (preprint 2002).


\bibitem {F}  {\sc Forstneri\v c,~F.} --- Extending proper holomorphic mappings of positive
codimension. {\em Invent. math.} {\bf 95}, no. 1, (1989), 31--61.



\bibitem{Hu94} {\sc Huang,~X.} --- On the mapping problem for algebraic real hypersurfaces in the complex spaces of different dimensions.
{\em Ann.~Inst.~Fourier $($Grenoble$)$} {\bf 44}, no. 2, (1994),
433--463.


\bibitem{HH} {\sc Huang,~X.} --- On some problems in several complex variables and CR geometry. {\em First International Congress of Chinese Mathematicians (Beijing, 1998)}, 383--396, {\em AMS/IP Stud. Adv. Math.} {\bf 20}, Amer. Math. Soc., Providence, RI, 2001.


\bibitem{la} {\sc Lamel, B.} --- Holomorphic maps of real submanifolds in complex spaces of different dimensions. {\em Pacific J. Math.} {\bf 201} (2001), no. 2, 357--387.

\bibitem{MMZ3}{\sc Meylan,~F.; Mir,~N.; Zaitsev,~D.} --- Approximation and convergence of formal CR-mappings.
{\em Internat. Math. Res. Notices} {\bf 4}, (2003), 211-242.

\bibitem{M00} {\sc Mir,~N.} --- Formal biholomorphic maps of real analytic hypersurfaces. {\em Math. Res. Lett.} {\bf 7} (2000), 343--359.

\bibitem{M02} {\sc Mir,~N.} --- On the convergence of formal mappings. {\em Comm. Anal. Geom.} {\bf 10}, no. 1, (2002), 23--59.

\bibitem{W} {\sc Webster,~S.M.} ---  On the mapping problem for
algebraic real hypersurfaces. {\em Invent. math.}, {\bf 43},
(1977), 53-68.

\end{thebibliography}
\end{document}